\documentclass{article}
\usepackage[utf8]{inputenc}
\usepackage{amsmath,amsthm,fullpage,graphicx}
\usepackage{hyperref}
\usepackage{listings}
\usepackage{authblk}

\title{Big Two and $n$-card poker probabilities}
\author[1]{Brian Wu}
\author[2]{Chai Wah Wu}
\affil[1]{New York University\thanks{briantehwu@gmail.com}}
\affil[2]{MIT-IBM Watson AI Lab, IBM T. J. Watson Research Center\thanks{cwwu@us.ibm.com}}

\date{August 15, 2023}

\begin{document}
\lstset{language=Python}

\maketitle

\begin{abstract}
    Between the poker hands of straight, flush, and full house, which hand is more common? In standard 5-card poker, the order from most common to least common is straight, flush, full house. The same order is true for 7-card poker such as Texas hold'em. However, is the same true for $n$-card poker for larger $n$? We study the probability of obtaining these various hands for $n$-card poker for various values of $n\geq 5$. In particular, we derive equations for the probability of flush, straight and full house and show that the probability of flush is less than a straight when $n \leq 11$, and is more than a straight when $n>11$.
    Similarly, we show that the probability of full house is less than a straight when $n \leq 19$, and is more than a straight when $n>19$.
    This means that for games such as Big Two where the ordering of $13$-card hands depends on the ordering in $5$-card poker, the ranking ordering does not follow the occurrence probability ordering, contrary to what intuition suggests.
\end{abstract}

\section{Introduction}
In traditional games of chance, value is typically correlated with the rarity of the outcome. In 5-card poker, the hands are ranked by their likelihood of occurrence. The royal flush, of which there are only four hands, is ranked the highest as it is the most uncommon. However, Texas hold'em, a popular variant of poker, does not follow this trend. While the same 5-card hands are constructed, players have seven cards to choose from. This changes the probabilities such that the lowest hand, the high card hand, is rarer than the one pair and two pair hands \cite{wiki:Poker_probability}. The ranking of the hands is the same as in 5-card poker even though the probabilities have changed. This quirk of the game is an exception to the rule: that value is correlated with rarity, and these exceptions are more common than one might think.

In the card game Big Two \cite{wiki:Big_two} a deck of 52 cards is distributed among the players (generally four). From here on, we will assume that Big Two is played with four players with each player being dealt 13 cards. Cards are played in groups (e.g. single cards, pairs, five card hands, etc). The value of a five card hand follows the corresponding value in 5-card poker, however, the five card hand must be straight or higher. While a straight is a more common hand in 5-card poker than a full house or a flush, is this still true for $13$-card poker? The 5-card poker hands of straight, flush and full house will be the focus of this paper as they are the most applicable to the card game Big Two.

The probability of obtaining specific hands in various variants of 5-card and 7-card poker are well studied \cite{wiki:Poker_probability}. For larger number of cards, such probabilities are less well known. Bill Butler \cite{durangobill:poker} calculated probabilities for $n$-card poker up to $n=10$ using a computer program implementing a brute-force approach to generate all possible hands. On the other hand, G. P. Michon provided a formula for the number of straight flush of $n$-card poker \cite{michon,oeisA143314}. 

However, general formulas for the number of flush, straight and full house are not available for arbitrary $n$.
The purpose of this paper is to give formulas and algorithms for computing the probability of flush, full house and straight in $n$-card poker and show that the relative order of these probabilities changes for different $n$.

\section{Notations and Preliminaries}
We assume a standard 52 card deck. The 4 suits are denoted as $s_i$ for $i=1,\cdots,4$. The number of $n$-card hands is $\binom{52}{n}$. For a multiset $R$, let $\nu(R)$ denote the number of permutations of the {\em nonzero} elements of $R$ given by a multinomial coefficient \cite{wiki:Permutation} and let $|R|$ to denote the number of {\em nonzero} elements of $R$ (counted with multiplicity). 

\section{Number of flush hands} \label{sec:flush}
A flush is defined as a set of cards with 5 or more cards of the same suit.
For now, we include into this set hands which have higher value such as four of a kind, straight flush and royal flush.
Let $r_i$ denote the number of cards of suit $s_i$. Let $S_n$
be the set of partitions of $n$ into a multiset of 4 integers in the range $[0,13]$. A hand is not a flush if and only if these integers are all less than $5$. 
For each multiset $S$ in $S_n$, a product of $3$ terms determines the total number of hands corresponding to $S$. The term $\binom{4}{|S|}$ describes the number of ways the multiset of nonzero elements of $S$ can be arranged among the 4 suits as $r_i$'s. The term $\nu(S)$ describes the number of arrangements of $S$, and finally $\Pi_i\binom{13}{r_i}$ describes the number of hands that has $S$ as the number of cards of the same suit.
Therefore, the number of non flush hands is equal to 
$\sum_{S\in S_n,\max(S)\leq 4}\binom{4}{|S|}\nu(S)\Pi_i\binom{13}{r_i}$.

As an illustrative example, consider the case of $n=7$ and $S=\{0,2,2,3\}$, i.e. the 7-card hands that have 2 cards of one suit, 2 cards of another suit and 3 cards of a third suit. Since $\max(S)=3 < 5$, these hands do not contain a flush. The nonzero elements (or equivalently the zero elements) of $S$ can be arranged in $\binom{4}{|S|}=4$ ways among the 4 suits, there are $\nu(S)=3$ permutations of the nonzero elements of $S$
and given an assignment of $S$ to $r_i$, there are $\binom{13}{2}\binom{13}{2}\binom{13}{3}$ ways to pick cards from these three suits making the number of hands corresponding to $S$ equal to $12\binom{13}{2}\binom{13}{2}\binom{13}{3}=20880288$.

The number of flush hands is then equal to 
$$\binom{52}{n}-\sum_{S=(r_1,r_2,r_3,r_4)\in S_n, r_i\leq 4}\binom{4}{|S|}\nu(\{r_1,r_2,r_3,r_4\})\binom{13}{r_1}\binom{13}{r_2}\binom{13}{r_3}\binom{13}{r_4}.$$
It is clear that by the pigeonhole principle for $n\geq 17$ there will be more than 5 cards with the same suit and the probability of a flush is $100\%$.
This remark shows that there are values of $n$ such that the probability of flush is larger than the probability of straight as it is easy to find 17 cards with no straight hands.

\section{Number of full house hands}
\label{sec:full}
A full house is defined as having 3 cards of the same rank and another 2 cards of the same rank. 
As before we include into this set hands which have higher value such as four of a kind.
We can use the same analysis and formulas in Section \ref{sec:flush} for flush hands
by simply swapping rank with suit. In particular, let $t_i$ denote the number of cards of rank $i$ where $1\leq i\leq 13$.
Let $H_n$
be the set of partitions of $n$ into 13 integers in the range
$[0,4]$.  Note that $|H_n|\leq |H_{26}|=104$. Let $T_n\subseteq H_n$ be the partitions in $H_n$ such that the largest term is larger than or equal to 3 and the second largest term is larger than or equal to $2$.
Then the number of full house hands is equal to 
$\sum_{T\in T_n}\binom{13}{|T|}\nu(T)\Pi_i\binom{4}{t_i}$.

Again, by the pigeonhole principle, we can always find a full house when $n\geq 27$, since there will be two suits with all the rank drawn along with a card from the third suit and the probability will be $100\%$. This remark shows that there are values of $n$ such that the probability of full house is larger than the probability of straight as it is easy to find 27 cards with no straight hands. This and the remark in the previous section show that when $n$ is small (e.g. $n=5$) straight is more common than full house or flush, but less common when $n$ is large. To get a clearer idea when and how often this switch occurs, we need to derive a formula for the number of straight hands. 

\section{Number of straight hands}
\label{sec:straight}
A straight is a list of 5 cards of sequential rank. We consider the standard high rules where the Ace can occur before Two and after a King card.
As before we include into this set hands which have higher value, such as four of a kind, as long as it contains a straight.
The first card of a straight can be an Ace, Two, $\cdots$ all the way to a Ten card.
Let $t_i$ denote the number of cards of rank $i$ and $H_n$ be as defined in Section \ref{sec:full}.
Let $g(k)\leq \binom{13}{k}$ be the number of combinations of $k$ distinct ranks (where each rank $r$ is in $1\leq r\leq 13$) which include a straight hand. 
Then, the number of straight hands is
$\sum_{H\in H_n}g(|H|)\nu(H)\Pi_i\binom{4}{t_i}$.

\subsection{Finding $g(k)$}

In order to prevent duplicate counting, let $W_i^k$ be the set of $k$ distinct ranks containing a straight where the smallest card in the straight has rank $i$. Therefore, $W_i^k\cap W_j^k=\emptyset$ if $i\neq j$.
It is clear that $g(k)=0$ for $k \leq 4$, since a straight hands needs $5$ distinct ranks. $g(5)=10$ since there are 10 ways that 5 cards of distinct ranks can be a straight hand, one for each of the 10 possible smallest cards. It is clear that $g(13)=1$ as there is only one way to have 13 distinct ranks and a straight hand is in this arrangement. Similarly $g(12) = 13$ since each of the $13$ arrangements of $12$ distinct ranks contains a straight hand. $g(11) = 77$ since there are $78$ arrangements of $11$ distinct ranks and only one such arrangement does not contains a straight (when the ranks do not contain a $5$ and a $10$ card).

Next, we show that $g(6)=71$.
For 6 distinct ranks, $|W_1^6|=8$, since the sixth card can be any card from rank 6 to rank 13. 
For $i>1$, the sixth rank can be any card except rank $i-1$, so $|W_i^6|=7$. This implies that $g(6) = 8+9\times 7=71$.
Similarly to compute $g(7)$, $|W_1^7| = \binom{8}{2}=28$, and for $i>1$, $|W_i^7|=\binom{7}{2}=21$, thus $g(7) = 28+9\times 21=217$, and
$g(8)=\binom{8}{3}+9\times\binom{7}{3}=56+9\times 35=371$.
To compute $g(9)$, we continue to use the formula $\binom{8}{4}+9\times\binom{7}{4}=70+9\times 35=385$, but need to remove cases of double counting. This occurs for $W_{10}^8$ which contains an Ace in the straight hand and picking the remaining $4$ cards to have rank $2,3,4,5$ will intersect with $W_1^8$. Thus $g(9) = 385-1=384$.

Finally, for $g(10)$, we use the formula $\binom{8}{5}+9\times\binom{7}{5}=56+9\times 21=245$ and list below 11 cases of double counting.
For $W_{10}^{10}$, the remaining cards being $(2,3,4,5,6)$, $(2,3,4,5,7)$, $(2,3,4,5,8)$, $(3,4,5,6,7)$, $(4,5,6,7,8)$ intersect with $W_1^{10}$, $W_2^{10}$, $W_3^{10}$ or $W_4^{10}$. 
For $W_9^{10}$, the remaining cards being $(1,2,3,4,5)$, $(2,3,4,5,6)$, $(3,4,5,6,7)$ intersect with $W_1^{10}$, $W_2^{10}$ and $W_3^{10}$ respectively.
For $W_8$, the remaining cards being $(1,2,3,4,5)$, $(2,3,4,5,6)$ intersect with $W_1^{10}$, and $W_2^{10}$ respectively.
For $W_7^{10}$, the remaining cards being $(1,2,3,4,5)$ intersects with $W_1^{10}$.
Therefore $g(10)=245-11=234$.

Thus we have shown that $g(k) = 0,0,0,0,10,71,217,371,384,234,77,13,1$ for $k=1,\cdots, 13$.

\subsection{A 11 card hand with no straight hand}
It is easy to see that a 11 card hand of the same suit does not contain a straight hand if and only if the ranks missing are 5 and 10. On the other hand, a 12 card hand of the same suit must contain a straight hand. 
When $n\geq 45$, by the generalized pigeonhole principle there is a suit with 12 or more cards and thus the number of straight hands is $\binom{52}{n}$ and the probability of a straight hand is $100\%$. Similarly,
when $n=44$, the only possibility that there is not a suit with 12 or more cards is when all 4 suits have exactly 11 cards and the only configuration where there is not a straight is if all the suits are missing rank 5 and 10. This means that the number of straight hands in $44$-card poker is $\binom{52}{44}-1=\binom{52}{8}-1$.

\section{Number of straight, flush and full house hands in $n$-card poker}
\begin{figure}[htbp]
\includegraphics[width=\textwidth]{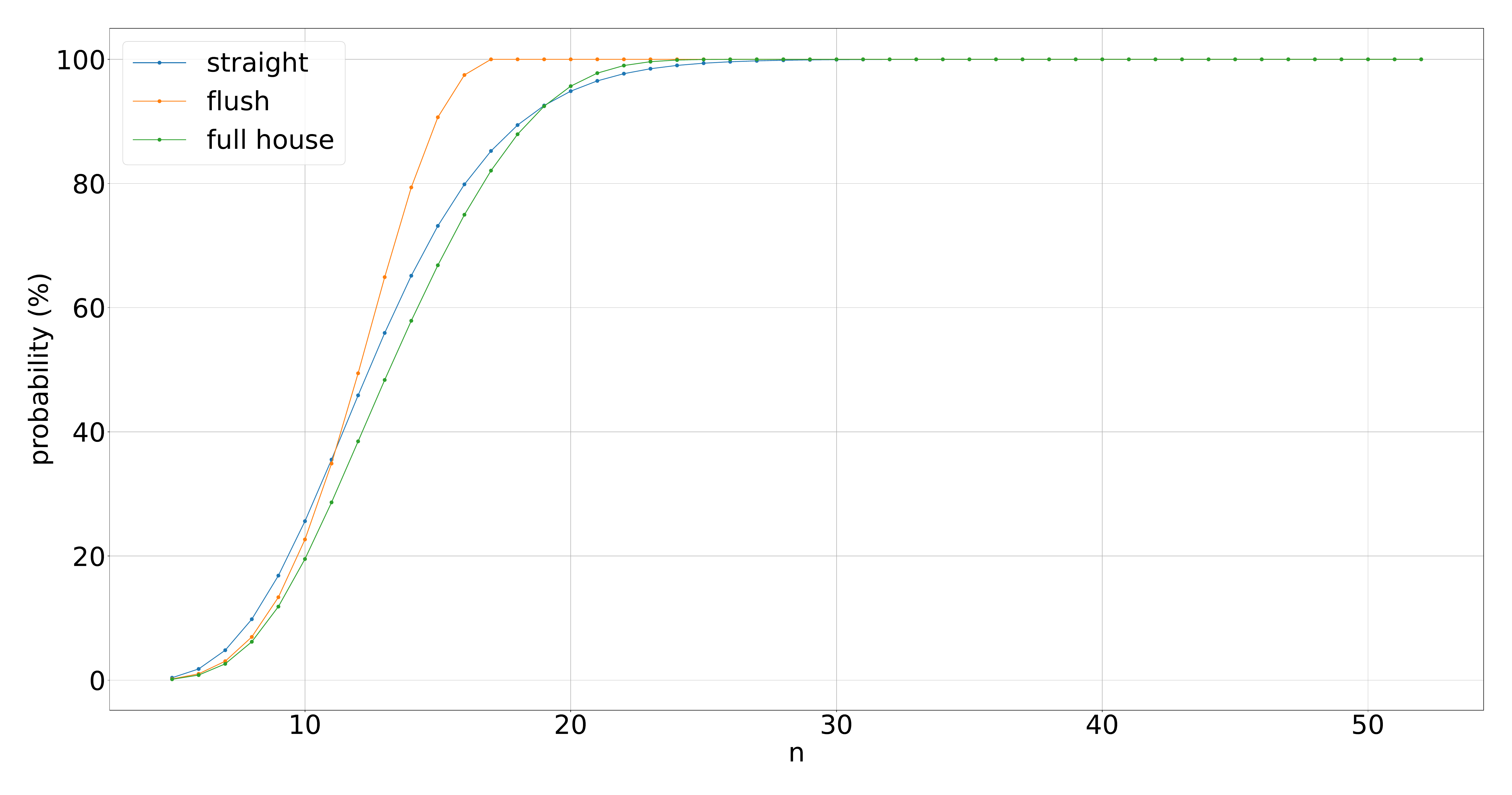}
\caption{Probability of various hands in $n$-card poker.}
\label{fig:poker}
\end{figure}

The values described in Sections \ref{sec:flush}-\ref{sec:straight} are computed and shown in Table \ref{tbl:poker}. 
Using the data in Table \ref{tbl:poker}, Figure \ref{fig:poker} shows the probability of occurrence of a straight, flush or full house in a random hand of $n$ cards. It shows that flush is always more probable than a full house. We see that for $n < 12$, straight hands are more common that full hands and flush hands. There is a crossover point at $n=12$, where for $n\geq 12$ flush hands are more common than straight hands. There is another crossover point at $n=20$, where for $n\geq 20$ full house hands are more common than straight hands.

\begin{table}[htbp]
\begin{center}
    \begin{tabular}{|c||c|c|c|c|}
    \hline
        $n$ & number of hands & number of straight & number of flush & number of full house\\ \hline\hline 
        5 & 2598960 & 10240 & 5148 & 3744 \\ \hline
        6 & 20358520 & 367616 & 207636 & 166920 \\ \hline
        7 & 133784560 & 6454272 & 4089228 & 3514992 \\ \hline
        8 & 752538150 & 73870336 & 52406640 & 46541430 \\ \hline
        9 & 3679075400 & 619588736 & 491448100 & 435926920 \\ \hline
        10 & 15820024220 & 4051217344 & 3585287134 & 3087272188 \\ \hline
        11 & 60403728840 & 21461806976 & 21076866408 & 17297489352 \\ \hline
        12 & 206379406870 & 94674009184 & 102014990714 & 79387982102 \\ \hline
        13 & 635013559600 & 355161047872 & 412247470340 & 307061893424 \\ \hline
        14 & 1768966344600 & 1152374363488 & 1404025311000 & 1024024781208 \\ \hline
        15 & 4481381406320 & 3279045142912 & 4063219805320 & 2994831165040 \\ \hline
        16 & 10363194502115 & 8276491135968 & 10101843501490 & 7769077277923 \\ \hline
        17 & 21945588357420 & 18706297925768 & 21945588357420 & 18011190562092 \\ \hline
        18 & 42671977361650 & 38154873848572 & 42671977361650 & 37522889445106 \\ \hline
        19 & 76360380541900 & 70680929691448 & 76360380541900 & 70598500404172 \\ \hline
        20 & 125994627894135 & 119535302593662 & 125994627894135 & 120551073059703 \\ \hline
        21 & 191991813933920 & 185328058520744 & 191991813933920 & 187726126771040 \\ \hline
        22 & 270533919634160 & 264282641858276 & 270533919634160 & 267830920323824 \\ \hline
        23 & 352870329957600 & 347526172985064 & 352870329957600 & 351537171709152 \\ \hline
        24 & 426384982032100 & 422213549653051 & 426384982032100 & 425903913135844 \\ \hline
        25 & 477551179875952 & 474573239602540 & 477551179875952 & 477437987194480 \\ \hline
        26 & 495918532948104 & 493971477605994 & 495918532948104 & 495905472254088 \\ \hline
        27 & 477551179875952 & 476384056580348 & 477551179875952 & 477551179875952 \\ \hline
        28 & 426384982032100 & 425743257091369 & 426384982032100 & 426384982032100 \\ \hline
        29 & 352870329957600 & 352546752515104 & 352870329957600 & 352870329957600 \\ \hline
        30 & 270533919634160 & 270384470692560 & 270533919634160 & 270533919634160 \\ \hline
        31 & 191991813933920 & 191928737926752 & 191991813933920 & 191991813933920 \\ \hline
        32 & 125994627894135 & 125970392012595 & 125994627894135 & 125994627894135 \\ \hline
        33 & 76360380541900 & 76351947773568 & 76360380541900 & 76360380541900 \\ \hline
        34 & 42671977361650 & 42669338815712 & 42671977361650 & 42671977361650 \\ \hline
        35 & 21945588357420 & 21944852453408 & 21945588357420 & 21945588357420 \\ \hline
        36 & 10363194502115 & 10363013522632 & 10363194502115 & 10363194502115 \\ \hline
        37 & 4481381406320 & 4481342680672 & 4481381406320 & 4481381406320 \\ \hline
        38 & 1768966344600 & 1768959253568 & 1768966344600 & 1768966344600 \\ \hline
        39 & 635013559600 & 635012471952 & 635013559600 & 635013559600 \\ \hline
        40 & 206379406870 & 206379271078 & 206379406870 & 206379406870 \\ \hline
        41 & 60403728840 & 60403715596 & 60403728840 & 60403728840 \\ \hline
        42 & 15820024220 & 15820023274 & 15820024220 & 15820024220 \\ \hline
        43 & 3679075400 & 3679075356 & 3679075400 & 3679075400 \\ \hline
        44 & 752538150 & 752538149 & 752538150 & 752538150 \\ \hline
        45 & 133784560 & 133784560 & 133784560 & 133784560 \\ \hline
        46 & 20358520 & 20358520 & 20358520 & 20358520 \\ \hline
        47 & 2598960 & 2598960 & 2598960 & 2598960 \\ \hline
        48 & 270725 & 270725 & 270725 & 270725 \\ \hline
        49 & 22100 & 22100 & 22100 & 22100 \\ \hline
        50 & 1326 & 1326 & 1326 & 1326 \\ \hline
        51 & 52 & 52 & 52 & 52 \\ \hline
        52 & 1 & 1 & 1 & 1 \\ \hline
\end{tabular}
\end{center}
\caption{Number of hands in $n$-card poker.}
\label{tbl:poker}
\end{table}

\subsection{Python code}
In Listings \ref{lst:straight}-\ref{lst:full} we provide the Python code used to generate the data in Table \ref{tbl:poker}

\begin{lstlisting}[label={lst:straight},frame=single,caption=Python code to compute the number of straight hands.]
# number of straight hands
from math import comb, prod, factorial
from sympy.utilities.iterables import partitions
def g(l):
    return (0,0,0,0,0,10,71,217,371,384,234,77,13,1)[l]
for n in range(5,53):
    k = 0
    for p in partitions(n,k=4,m=13):
        ps = []
        for d in p:
            ps.extend([d]*p[d])
        q = sorted(ps)
        k += g(len(ps)) \
            *(factorial(len(ps))//prod(factorial(q) for q in p.values())) \
            *prod(comb(4,q) for q in ps)
    print(n,k)
\end{lstlisting}

\begin{lstlisting}[label={lst:flush},frame=single,caption=Python code to compute the number of flush hands.]
# number of flush hands
from math import comb, prod, factorial
from sympy.utilities.iterables import partitions
for n in range(5,53):
    k = 0
    for p in partitions(n,k=13,m=4):
        ps = []
        for d in p:
            ps.extend([d]*p[d])
        q = sorted(ps)
        if (max(q)<=4):
            k += comb(4,len(ps)) \
                *(factorial(len(ps))//prod(factorial(q) for q in p.values())) \
                *prod(comb(13,q) for q in ps)
    print(n,comb(52,n)-k)
\end{lstlisting}

\begin{lstlisting}[label={lst:full},frame=single,caption=Python code to compute the number of full house hands.]
# number of full house hands
from math import comb, prod, factorial
from sympy.utilities.iterables import partitions
for n in range(5,53):
    k = 0
    for p in partitions(n,k=4,m=13):
        ps = []
        for d in p:
            ps.extend([d]*p[d])
        q = sorted(ps)
        if (q[-1] >= 3 and q[-2] >= 2):
            k += comb(13,len(ps)) \
                *(factorial(len(ps))//prod(factorial(q) for q in p.values())) \
                *prod(comb(4,q) for q in ps)
    print(n,k)
\end{lstlisting}

\section{Conclusions}
We considered the poker hands of straight, flush, and full house in $n$-card poker. 
We show that a straight hand can be more or less probable than a flush or a full house depending on the value of $n$. This shows that in four player Big Two, the odds of drawing a hand with a flush are higher than the odds of drawing a hand with a straight, contrary to how they are scored.

\bibliographystyle{ieeetr}
\bibliography{poker}

\end{document}